\documentclass[10pt,reqno]{amsart}
\usepackage{amsfonts,amsmath,amssymb,amsbsy,amstext,amsthm}
\usepackage{accents,color,enumerate,enumitem,float,fullpage,verbatim}

\usepackage{url}
\usepackage{microtype}
\usepackage[colorlinks=true,hyperindex, linkcolor=magenta, pagebackref=false, citecolor=cyan]{hyperref}
\usepackage[alphabetic,lite,backrefs]{amsrefs}

\usepackage{eucal,bm,kpfonts,mathbbol, parskip}

\usepackage{tikz,tikz-cd}
\usetikzlibrary{positioning, matrix, shapes}
\usetikzlibrary{decorations.pathmorphing}

\usepackage{lscape}

\usepackage{titlesec}
\titleformat{\section}[block]{\large\scshape\bfseries\filcenter}{\thesection.}{1em}{}		
\titleformat{\subsection}[hang]{\large\scshape\bfseries}{\thesubsection}{1em}{}			
\titleformat{\subsubsection}[hang]{\large\scshape\bfseries}{\thesubsubsection}{1em}{}		

\usepackage{multirow}
\usepackage{array}
\usepackage{booktabs}
\newcolumntype{M}[1]{>{\centering\arraybackslash}m{#1}}
\newcolumntype{N}{@{}m{0pt}@{}}
\usepackage{diagbox}
\usepackage{cancel}

\newtheorem{lemma}{Lemma}[section]
\newtheorem{theorem}[lemma]{Theorem}

\newtheorem{definition}[lemma]{Definition}

\theoremstyle{remark}
\newtheorem{remark}[lemma]{Remark}
\newtheorem{example}[lemma]{Example}

\newcommand{\Pic}{\operatorname{Pic}}

\newcommand{\Cox}{\operatorname{Cox}}

\newcommand{\rank}{\operatorname{rank}}

\newcommand{\doot}{\bullet}

\renewcommand{\aa}{\mathbf a}

\newcommand{\dd}{\mathbf d}

\newcommand{\nn}{\mathbf n}

\newcommand{\fF}{\mathbf F}

\renewcommand{\O}{\mathcal{O}}

\newcommand{\NN}{\mathbb{N}}
\newcommand{\PP}{\mathbb{P}}

\newcommand{\ZZ}{\mathbb{Z}}

\title{The Virtual Resolutions Package for Macaulay2}

\author[Ayah Almousa]{Ayah Almousa}
\address{Department of Mathematics, Cornell University, Ithaca, New York, 14853}
\email{\href{mailto:aka66@cornell.edu}{aka66@cornell.edu}}
\urladdr{\url{http://math.cornell.edu/~aalmousa}}

\author[Juliette Bruce]{Juliette Bruce}
\address{Department of Mathematics, University of Wisconsin, Madison, Wisconsin, 53706}
\email{\href{mailto:juliette.bruce@math.wisc.edu}{juliette.bruce@math.wisc.edu}}
\urladdr{\url{http://math.wisc.edu/~juliettebruce/}}

\author[M.C. Loper]{Michael C. Loper}
\address{School of Mathematics, University of Minnesota, Minneapolis, Minnesota, 55455}
\email{\href{mailto:loper012@umn.edu}{loper012@umn.edu}}
\urladdr{\url{http://math.umn.edu/~loper012/}}

\author{Mahrud Sayrafi}
\address{School of Mathematics, University of Minnesota, Minneapolis, Minnesota, 55455}
\email{\href{mailto:mahrud@umn.edu}{mahrud@umn.edu}}
\urladdr{\url{http://math.umn.edu/~mahrud/}}


\begin{document} 
\thanks{AA was partially supported by the NSF GRFP under Grant No. DGE-1650441.}
\thanks{JB was partially supported by the NSF GRFP under Grant No. DGE-1256259 and NSF grant DMS-1502553.}
\thanks{ML was partially supported by the NSF RTG grant DMS-1745638}

\maketitle
\vspace{-1cm}
\begin{abstract}
We introduce the \texttt{VirtualResolutions} package for the computer algebra system \textit{Macaulay2}. This package has tools to construct, display, and study virtual resolutions for products of projective spaces. The package also has tools for generating curves in $\PP^1\times\PP^2$, providing sources of interesting virtual resolutions.
\end{abstract}

\setcounter{section}{1}


Recently, Berkesch, Erman, and Smith introduced the notion of virtual resolutions for subvarieties of smooth projective toric varieties as an analogue to minimal graded free resolutions for subvarieties of projective space \cite{berkesch17}. A virtual resolution is a complex of finitely generated graded free modules over the Cox ring of a smooth projective toric variety that becomes exact upon passing to the corresponding complex of coherent sheaves.

While graded minimal free resolutions are useful for studying quasicoherent sheaves on projective spaces, when working over a product of projective spaces or, more generally, over smooth projective toric varieties, they are often long and cumbersome to compute.
By allowing a limited amount of homology, virtual resolutions offer a more flexible alternative for studying toric subvarieties when compared to graded minimal free resolutions. This article introduces the \texttt{VirtualResolutions} package for \textit{Macaulay2} \cite{M2}, which provides tools for constructing and studying virtual resolutions over products of projective spaces. 

The \texttt{VirtualResolutions} project began in 2018 at the \textit{Macaulay2 Workshop at the University of Wisconsin - Madison}, building on previous work of Christine Berkesch, David Eisenbud, Daniel Erman, and Gregory G. Smith. Along with them, the authors of this paper developed the \texttt{VirtualResolutions} package to generate examples of virtual resolutions in products of projective spaces. In particular, the package provides tools to generate examples of virtual resolutions and to check whether chain complexes are virtual resolutions. These methods are introduced and demonstrated in Section 2.

In Section 3, we concentrate on virtual resolutions arising from curves in $\PP^1 \times \PP^2$. This is the simplest case after points in a product of projective spaces. We demonstrate several methods for finding the defining ideals of curves in the Cox ring of $\PP^1 \times \PP^2$, including monomial curves, rational curves, and curves from $\PP^3$.

\section*{Acknowledgments}  

We would like to thank Christine Berkesch, David Eisenbud, Daniel Erman, Gregory G. Smith, and Mike Stillman for helping to develop the \texttt{VirtualResolutions} package. Further, we especially thank Christine Berkesch for reading a draft of this article and providing many helpful comments. We would also like to thank the referees for many helpful comments.

Work on this article and the associated \textit{Macaulay2} package began at the \textit{Macaulay2 Workshop at the University of Wisconsin - Madison} funded by NSF grant DMS-1812462. Additional progress was made during the \textit{Graduate Workshop in Commutative Algebra for Women and Mathematicians of Minority Genders} funded by NSF grant DMS-1908799. The final portion of this article and package was written during a Coding Sprint hosted by the Institute for Mathematics and its Applications (IMA), and we are grateful for their hospitality.

\section{Constructing Virtual Resolutions}\label{sec:vres}

If $X$ is a smooth projective toric variety, we denote its $\Pic(X)$-graded Cox ring, as defined in \cite{cox95}, by $\Cox(X)$ (see also Section~5.2 of \cite{coxLittleSchenck11}). Denote the associated irrelevant ideal of $X$ by $B\subset \Cox(X)$. A virtual resolution of a $\Pic(X)$-graded module over $\Cox(X)$ is defined as follows:
 

\begin{definition}\cite{berkesch17}*{Definition~1.1}
 A \textbf{virtual resolution} of a $\Pic(X)$-graded $\Cox(X)$-module $M$ is a chain complex of $\Pic(X)$-graded free $\Cox(X)$-modules
\begin{center}
\begin{tikzcd}[column sep = 3em,ampersand replacement=\&]
\fF_{\doot}\coloneqq\big[ F_{0} \& \lar F_{1} \& \lar F_{2} \& \lar \cdots \big]
\end{tikzcd}
 \end{center}
 such that the corresponding complex of $\O_{X}$-modules is a locally free resolution of the sheaf $\widetilde{M}$. 
\end{definition}

One can rephrase the definition of a virtual resolution in a way that is more practical for computations: a complex $\fF_{\doot}$ of $\Pic(X)$-graded free $\Cox(X)$-modules is a virtual resolution of $M$ if and only if the following two conditions are satisfied:
\begin{enumerate}
    \item The $B$-saturation of $H_0(\fF_{\doot})$ is isomorphic to the $B$-saturation of $M$. 
    \item The homology module $H_{i}(\fF_{\doot})$ is supported only on $B$ for all $i>0$.
\end{enumerate}

For the remainder of this paper, we will focus on the case when $X$ is a product of projective spaces, since \cite{berkesch17} demonstrates several ways to construct virtual resolutions in this case. We let $\nn = (n_1,n_2,...,n_r)\in (\ZZ_{\geq0})^r$ and $\PP^{\nn} = \PP^{n_1} \times \PP^{n_2} \times \cdots \times \PP^{n_r}$. Further, we write $S$ for $\Cox(\PP^{\nn})$, which is graded by $\Pic(\PP^{\nn})\cong\ZZ^r$.

Given this setup, Berkesch, Erman, and Smith show that it is possible to construct a virtual resolution of a $\ZZ^{r}$-graded $S$-module $M$ from the graded minimal free resolution of $M$ and $\dd\in \ZZ^r$ an element of the multigraded Castelnuovo-Mumford regularity as introduced in \cite{maclaganSmith04}. They call a virtual resolution constructed in this manner the \textbf{virtual resolution of the pair} $(M,\dd)$.

\begin{theorem}\cite{berkesch17}*{Theorem~1.3}\label{thm:virtualofpair} Let $M$ be a finitely-generated $\ZZ^r$-graded $B$-saturated $S$-module that is $\dd$-regular. If $G$ is the free subcomplex of a minimal free resolution of $M$ consisting of all summands generated in degrees at most $\dd+\nn$, then $G$ is a virtual resolution of $M$.

\end{theorem}


Notice that the virtual resolution of the pair $(M,\dd)$ is a subcomplex of the minimal free resolution, so its length is at most the projective dimension of $M$. Also note that we are justified in saying ``the'' virtual resolution of the pair $(M,\dd)$ because this complex is unique up to isomorphism.

\begin{example}\label{virtualOfPair}
  \label{ex:3pts}
  Consider three points $\left([1:1],[1:4]\right)$, $\left([1:2],[1:5]\right)$, and $\left([1:3],[1:6]\right)$ in $\PP^1 \times \PP^1$. Using \textit{Macaulay2}, we may compute the $B$-saturated ideal defining these three points.
\begin{verbatim}
    i1 : needsPackage "VirtualResolutions";
    i2 : X = toricProjectiveSpace(1) ** toricProjectiveSpace(1);
    i3 : S = ring X; B = ideal X;
    o4 : Ideal of S
    i5 : J = saturate(intersect(
             ideal(x_1 - x_0, x_3 - 4*x_2),
             ideal(x_1 - 2*x_0, x_3 - 5*x_2),
            ideal(x_1 - 3*x_0, x_3 - 6*x_2)), B);
    o5 : Ideal of S
\end{verbatim}

One can show the $\ZZ^2$-graded minimal free resolution of $S/J$ is the following complex.
\begin{center}
\begin{tikzcd}[column sep = 2.5em]
S & \lar \begin{matrix} S(0,-3) \\ \oplus \\ S(-1,-2) \\ \oplus \\ S(-2,-1) \\ \oplus \\ S(-3,0) \\ \oplus \\ S(-1,-1) \end{matrix}
  & \lar \begin{matrix} S(-1,-3)^2 \\ \oplus \\ S(-2,-2)^2 \\ \oplus \\ S(-3,-1)^2 \end{matrix}
  & \lar \begin{matrix} S(-2,-3) \\ \oplus \\ S(-3,-2) \end{matrix}
  & \lar 0
\end{tikzcd}
\end{center}

We may view the multigraded Betti table for the above graded minimal free resolution.

\begin{verbatim}
    i6 : minres = res J;
    i7 : multigraded betti minres
            0             1               2         3
    o7 = 0: 1             .               .         .
         2: .            ab               .         .
         3: . a3+a2b+ab2+b3               .         .
         4: .             . 2a3b+2a2b2+2ab3         .
         5: .             .               . a3b2+a2b3
    o7 : MultigradedBettiTally
\end{verbatim}

In order to compute a virtual resolution of $S/J$, we find an element of the multigraded Castelnuovo-Mumford regularity of the module \cite{maclaganSmith04}*{Definition~4.1}. The minimal elements in the multigraded regularity of $S/J$ can be computed via the \texttt{multigradedRegularity} command.

\begin{verbatim}
    i8 : multigradedRegularity(X, S^1/J)
    o8 = {{0, 2}, {1, 1}, {2, 0}}
    o8 : List
\end{verbatim}


This computation relies on the fact that a $B$-saturated $S$-module $M$ over a product of projective spaces $\PP^\nn$ is $\dd$-regular provided that the following conditions are satisfied: 
\begin{enumerate}
\item the Hilbert function $H(M, \dd)$ agrees with the Hilbert polynomial $P_M(\dd)$; 
\item $H^i\left(\PP^\nn, \widetilde{M}\left(\aa\right)\right) = 0$ for all $i>0$ and twists $\aa$ such that $d_j-i\leq a_j$. 
\end{enumerate}
In particular, the function \texttt{multigradedRegularity} utilizes the \textit{Macaulay2} package \texttt{TateOnProducts}  \cite{eisenbudErmanSchreyer19} to compute the sheaf cohomology of twists of $\widetilde{M}$, which allows one to find elements of the multigraded regularity of the sheaf \cite{eisenbudErmanSchreyer15}*{Proposition~3.11}.


The function \texttt{virtualOfPair} implements Theorem \ref{thm:virtualofpair} to compute the virtual resolution of a pair. We call the function below in order to compute the virtual resolution of the pair $(S/J,(2,0))$. Note that since we must remove all twists generated in degrees not less than or equal to $\nn +\dd$ from the graded minimal free resolution, we input the element $(3,1)=(1,1)+(2,0)$ in  \texttt{virtualOfPair}.

\begin{verbatim}
    i9 : vres = virtualOfPair(res J, {{3, 1}})

          1      3      2
    o9 = S  <-- S  <-- S  <-- 0
                          
         0      1      2      3 
\end{verbatim}
\pagebreak
\begin{verbatim}
    i10 : multigraded betti vres

             0      1    2
    o10 = 0: 1      .    .
          2: .     ab    .
          3: . a3+a2b    .
          4: .      . 2a3b
    o10 : MultigradedBettiTally
\end{verbatim}

The above virtual resolution of the pair $(S/J, (2,0))$ is shorter and thinner than the graded minimal free resolution of $S/J$. 

When a minimal free resolution is already known or provided as input, \texttt{virtualOfPair} takes the appropriate subcomplex of that resolution. Otherwise, \texttt{virtualOfPair} computes a virtual resolution by using the induced Schreyer orders to iteratively compute syzygies in the desired degrees. Note that while the resulting virtual resolutions from the two strategies may have different differentials, the chain complexes are unique up to isomorphism.

Continuing the example of three points in $\mathbb{P}^1\times\mathbb{P}^1$, we again use \texttt{virtualOfPair} to compute the virtual resolution of $S/J$. In this case, \textit{Macaulay2} does not have the minimal free resolution cached, so \texttt{virtualOfPair} uses Schreyer's method to obtain a virtual resolution isomorphic to the one above.

\begin{verbatim}
    i11 : multigraded betti virtualOfPair(S^1/J, {{3, 1}})

             0      1    2
    o11 = 0: 1      .    .
          2: .     ab    .
          3: . a3+a2b    .
          4: .      . 2a3b
    o11 : MultigradedBettiTally
\end{verbatim}

\begin{remark}
  In larger experiments, Schreyer's method is significantly more efficient than first computing a minimal free resolution. For instance, finding a virtual resolution for the saturated ideal of $5$ points in $(\PP^1)^4$ using Schreyer's method is two orders of magnitude faster. These time savings grow even more dramatically as the number of points or number of factors of $\PP^1$ are increased.
\end{remark}

We may check that \texttt{vres} is indeed a virtual resolution by using the \texttt{isVirtual} function.

\begin{verbatim}
    i12 : isVirtual(B, vres)
    o12 = true
\end{verbatim}

By default, \texttt{isVirtual} checks whether the homology of the given chain complex is supported only on the irrelevant ideal. More specifically, \texttt{isVirtual} checks whether the annihilator of the homology of the given chain complex saturates to the entire ring.

\begin{remark}
Computing the saturation of an ideal is a critical but computationally costly step in many aspects of studying virtual resolutions. As such, our package utilizes new saturation methods which are under development by Michael Stillman for future release in \textit{Macaulay2}. These methods are currently stored in the auxiliary file \texttt{Colon.m2}, but will be removed once they have been released separately.
\end{remark}

As an alternative approach, we also implement the method presented in Theorem 1.3 of \cite{loper19} for checking whether a complex is virtual. This is done by setting \verb+Strategy => Determinantal+. In this case, two conditions are checked: one involving the ranks of the maps in the chain complex, and the other involving the depths of the $B$-saturated ideals of minors of the maps. Typically the default strategy is faster than the determinantal strategy, as the determinantal strategy  must compute the ideal of minors for each map in the chain complex, which is generally computationally expensive.
\end{example}

Another way to generate virtual resolutions is by using the following theorem of \cite{berkesch17}, which provides a method for producing virtual resolutions of ideals defining zero-dimensional subschemes of $\PP^{\nn}$.

\begin{theorem}\cite{berkesch17}*{Theorem~4.1}\label{thm:viafatpoint} If $Z\subset \PP^\nn$ is a zero-dimensional scheme and $I$ is the corresponding $B$-saturated $S$-ideal, then there exists $\aa\in\NN^r$ with $a_r = 0$ such that the minimal free resolution of $S/(I\cap B^\aa)$ has length equal to $|\nn| = \dim\PP^\nn$. Moreover, any $\aa\in\NN^r$ with $a_r = 0$ and other entries sufficiently positive yields such a virtual resolution of $S/I$.
\end{theorem}

The function \texttt{resolveViaFatPoint} computes $B^{\aa} = \bigcap_{i = 1} ^r (x_{i,0},x_{i,1},\ldots,x_{i,n_i})^{a_i}$, intersects it with the input ideal $J$, and computes the minimal free resolution of $S/(J\cap B^{\aa})$.

\begin{example} \label{resolveViaFatPoint}
Again consider the ideal of three points in $\PP^1\times \PP^1$ as in Example \ref{virtualOfPair}. We use the function \texttt{resolveViaFatPoint} to compute a virtual resolution of $J$.

\begin{verbatim}
    i13 : C = resolveViaFatPoint(J, B, {2, 0})

           1      4      3
    o13 = S  <-- S  <-- S  <-- 0

          0      1      2      3
      
    o13 : ChainComplex
    i14 : multigraded betti C

             0       1    2
    o14 = 0: 1       .    .
          3: . a3+3a2b    .
          4: .       . 3a3b

    o14 : MultigradedBettiTally
    i15 : isVirtual(B, C)
    o15 = true
\end{verbatim}

Note that the virtual resolutions obtained by
\texttt{resolveViaFatPoint} may not be a subcomplex of the minimal resolution of $M$ (as demonstrated by Examples~\ref{virtualOfPair} and~\ref{resolveViaFatPoint}), and thus need not occur as the virtual resolution of
$(M, \dd)$ for a given regularity $\dd$. Thus, \texttt{resolveViaFatPoint} is a useful method for producing new and interesting examples of virtual resolutions.
\end{example}

\section{Constructing Curves in \texorpdfstring{$\PP^1 \times \PP^2$}{P1xP2}}\label{sec:curves}
One source of potentially interesting virtual resolutions comes from studying curves in $\PP^1\times\PP^2$. With this in mind, the \texttt{VirtualResolutions} package contains functions for constructing a limited class of curves in $\PP^1\times\PP^2$, which we believe may be of interest to other researchers. In particular, there are functions for generating monomial curves, rational curves, and curves arising from curves in $\PP^3$.

The main function in this direction is \texttt{curveFromP3toP1P2}. Given the defining ideal of a curve in $\PP^3$, this function returns the defining ideal of a curve in $\PP^1\times\PP^2$ constructed in the following way: given projections $\pi_{1}:\PP^3\dashrightarrow\PP^1$ and $\pi_{2}:\PP^3\dashrightarrow\PP^2$, there is an induced rational map $\psi:\PP^3\dashrightarrow \PP^1\times\PP^2$, and this is the map under which we are computing the image of our curve $C$. 

The projections that \texttt{curveFromP3toP1P2} uses are fixed. In particular, it uses the coordinate projections 
\[
\pi_1\left([z_0:z_1:z_2:z_3]\right)=[z_0:z_1] \quad \text{and}  \quad \pi_{2}\left([z_0:z_1:z_2:z_3]\right)=[z_1:z_2:z_3].
\]
As one might wish to preserve the degree of the curve $C\subset\PP^1\times\PP^2$, the function \texttt{curveFromP3toP1P2} has an option called \texttt{PreserveDegree}. When this option is set to \texttt{true}, the function will return an error if the given curve intersects the base locus of these projections. 

\begin{example}\label{ex:virtualCurve}
The code below uses \texttt{curveFromP3toP1P2} to construct a curve in $\PP^1\times\PP^2$ from the defining ideal of the twisted cubic $C\subset \PP^3$. 
\begin{verbatim}
    i16 : R = ZZ/101[z_0, z_1, z_2, z_3];
    i17 : I = ideal(z_0*z_2-z_1^2, z_1*z_3-z_2^2, z_0*z_3-z_1*z_2);
    o17 : Ideal of R
    i18 : J = curveFromP3toP1P2 I
                  2
    o18 = ideal (x    - x   x   , - x   x    + x   x   , - x   x    + x   x   )
                  1,1    1,0 1,2     0,1 1,1    0,0 1,2     0,1 1,0    0,0 1,1

                   ZZ
    o18 : Ideal of ---[x   , x   , x   , x   , x   ]
                   101  0,0   0,1   1,0   1,1   1,2
\end{verbatim}
We can check that the ideal $J$ defines a curve in $\PP^1\times\PP^2$ by computing its dimension. Note that since $J$ is an ideal in the Cox ring of $\PP^1\times\PP^2$, we expect $J$ to be three dimensional, since in general the dimension of a subscheme defined by a $\Pic(X)$-graded ideal $J\subset \Cox(X)$ is given by $\dim J-\rank \Pic(X)$.
\begin{verbatim}
    i19 : dim J
    o19 = 3
\end{verbatim}
\end{example}

The function \texttt{randomCurveP1P2} constructs a random curve in $\PP^1\times \PP^2$ by calling the \texttt{curve} function from the package \texttt{SpaceCurves} \cites{zhang18} to randomly generate a curve $C\subset\PP^3$ of specified degree and genus; see \cite{zhang18}*{Section~2} for a detailed discussion of the inner workings of \texttt{curve}. The function \texttt{randomCurveP1P2} then constructs a curve $C'\subset \PP^1\times\PP^2$ from $C$ using
\texttt{curveFromP3toP1P2}.


\begin{example}
Using the \texttt{randomCurveP1P2} function we produce a random curve in $\PP^1\times\PP^2$ coming from a curve of degree $7$ and genus $3$ in $\PP^3$. Additionally, we check that the resulting ideal defines a curve.
\begin{verbatim}
    i20 : I = randomCurveP1P2(7, 3);
                    ZZ
    o20 : Ideal of ---[x   , x   , x   , x   , x   ]
                   101  0,0   0,1   1,0   1,1   1,2
    i21 : S = ring I;
    i22 : dim I
    o22 = 3
\end{verbatim}
Using the \texttt{multigradedRegularity} function, one can see that $(2,3)$ is a minimal element in the multigraded regularity of this curve. Note that we first saturate the ideal defining our curve; this is necessary to ensure that the output of \texttt{multigradedRegularity} is correct.
\begin{verbatim}
    i23 : B = intersect(ideal(x_(0,0), x_(0,1)), ideal(x_(1,0), x_(1,1), x_(1,2)));
    o23 : Ideal of S
    i24 : J = saturate(I, B);
    o24 : Ideal of S
    i25 : multigradedRegularity(S, S^1/J)
    o25 = {{1, 4}, {2, 3}}
    o25 : List
\end{verbatim}
Finally, we can use this element of the multigraded regularity to compute a virtual resolution of $S/J$. \pagebreak
\begin{verbatim}
    i26 : minres = res J;
    i27 : vres = virtualOfPair(J, {{3, 5}})

           1      11      18      8
    o27 = S  <-- S   <-- S   <-- S  <-- 0
                                     
          0      1       2       3      4

    o27 : ChainComplex
\end{verbatim}
Comparing the multigraded Betti tables of these two resolutions we see that once again the virtual resolution is shorter and less complicated than the graded minimal free resolution. 
\begin{verbatim}
    i28 : multigraded betti minres
             0               1                       2                 3          4
    o28 = 0: 1               .                       .                 .          .
          5: .      a3b2+3a2b3                       .                 .          .
          6: . a4b2+2a2b4+5ab5                       .                 .          .
          7: .              b7 3a4b3+6a3b4+12a2b5+3ab6                 .          .
          8: .               .                    2ab7 3a4b4+8a3b5+6a2b6          .
          9: .               .                       .              a2b7 a4b5+3a3b6
    o28 : MultigradedBettiTally
    i29 : multigraded betti vres
             0          1            2     3
    o29 = 0: 1          .            .     .
          5: . a3b2+3a2b3            .     .
          6: . 2a2b4+5ab5            .     .
          7: .          . 6a3b4+12a2b5     .
          8: .          .            . 8a3b5
    o29 : MultigradedBettiTally
\end{verbatim}
\end{example}

\begin{remark}\label{rem:randomReg}
Since the curve generated by \texttt{randomCurveP1P2} is, in some sense, random, one will often get different virtual resolutions when running the above example. As an example of this phenomena, working over a slightly larger finite field than in Example~\ref{ex:virtualCurve}, we compute the multigraded regularity for 500 curves arising from curves of genus 3 and degree 7 in $\PP^3$. The resulting distribution of multigraded regularities is shown below.
\begin{verbatim}
    i30 : tally apply(500,i->(
              I := randomCurveP1P2(7, 3, ZZ/32003);
              S := ring I;
              B := intersect(ideal(x_(0,0), x_(0,1)),
              ideal(x_(1,0), x_(1,1), x_(1,2)));
              J := saturate(I, B);
              multigradedRegularity(S, S^1/J)
           ))
    o30 = Tally{{{1, 4}, {2, 3}} => 167}
                 {{1, 5}} => 170
                 {{1, 6}} => 163
    o30 : Tally
\end{verbatim}
\end{remark}

\begin{remark}
It would be interesting to have an understanding of how the geometry of the curves in $\PP^3$ affects multigraded regularity of the resulting curves in $\PP^1\times\PP^2$. For example, in the computation in Remark~\ref{rem:randomReg} every smooth curve of genus 3 and degree 7 lies on either a smooth cubic surface or a rational quartic surface with a double line. In fact, using the notation from \cite{zhang18}, each curve of genus 3 and degree 7 lies in one of three possible divisor classes: $(4, 1, 1, 1, 1, 1, 0)$, $(5, 3, 1, 1, 1, 1, 1)$, and $(6, 1, 2, 2, 2, 2, 2, 2, 2, 1)$, the first two being divisors on a smooth cubic surface and the last being a divisor on a rational quartic surface with a double line. An experiment similar to the one in Remark~\ref{rem:randomReg} suggests that the following relationship between divisors and multigraded regularity:
\begin{enumerate}
    \item The minimal elements in the multigraded regularity of any curve in $\PP^1\times\PP^2$ arising from a smooth curve in $\PP^3$ lying on smooth cubic surface in the divisor class $(4, 1, 1, 1, 1, 1, 0)$ are $\{(1,6)\}$.
    \item The minimal elements in the multigraded regularity of any curve in $\PP^1\times\PP^2$ arising from a smooth curve in $\PP^3$ lying on smooth cubic surface in the divisor class $(5, 3, 1, 1, 1, 1, 1)$ are $\{(1,5)\}$.
    \item The minimal elements in the multigraded regularity of any curve in $\PP^1\times\PP^2$ arising from a smooth curve in $\PP^3$ lying on a rational quartic surfaces with a double line in the divisor class $(6, 1, 2, 2, 2, 2, 2, 2, 2, 1)$ are $\{(1,4), (2,3)\}$.
\end{enumerate}
How this generalizes to all curves though remains unclear, but we hope these functions will stimulate others to consider this and other similar problems. 
\end{remark}

The functions \texttt{randomRationalCurve} and \texttt{randomMonomialCurve} are alternative methods for generating curves in $\PP^1\times\PP^2$. Given two positive integers $d$ and $e$, these functions construct the defining ideal of the curve arising as the image of the map
\begin{center}
\begin{tikzcd}[column sep = 3em]
\PP^1 \rar & \PP^1\times\PP^2 & &\text{given by} & & \left[t_0 : t_1\right] \rar[mapsto] & \left(\left[f_0:f_1\right],\left[g_0:g_1:g_2\right]\right)
\end{tikzcd}
\end{center}
where the $f_i$ and $g_i$ are forms in $\mathbb{K}[t_0,t_1]$ of degrees $d$ and $e$, respectively. For \texttt{randomMonomialCurve}, these forms are chosen to be monomials, while for \texttt{randomRationalCurve} they can be any possible forms of the correct degree.

\begin{example}
Here we use the \texttt{randomRationalCurve} function to construct a rational curve of bidegree $(5,7)$. Once again we verify that it is, in fact, a curve by computing its dimension. 
\begin{verbatim}
    i31 : I = randomRationalCurve(5, 7);
                    ZZ
    o31 : Ideal of ---[x , x , y , y , y ]
                   101  0   1   0   1   2
    i32 : dim I 
    o32 = 3
\end{verbatim}
\end{example}

\end{document}